\documentclass[a4paper,10pt]{article}

\usepackage[utf8]{inputenc}

\usepackage[english ]{babel}
 
\usepackage{xcolor}
\usepackage{bbding}

\usepackage[utf8]{inputenc}

\usepackage[all]{xy}

\usepackage{tikz-cd}
\usepackage[all]{xy}
\usepackage{amsfonts}
\usepackage{amssymb}
\usepackage{graphicx}
\usepackage{mathtools}
\usepackage{skak} 
\usepackage{foekfont}

\usepackage{calligra}
\usepackage[T1]{fontenc}

\usepackage{amsmath,mathdots}

\usepackage[utf8]{inputenc}
\usepackage{devanagari}

\usepackage{mathtools}
\usepackage{mathdots}
\usepackage{tikz}

\definecolor{amber(sae/ece)}{rgb}{1.0, 0.49, 0.0}
\newfont{\rsfsten}{rsfs10 scaled 1200}
\usepackage{MnSymbol}
\usepackage{tikz}
\usepackage{amscd}
\usepackage{MnSymbol}
\usepackage{wasysym}

\usepackage{mathrsfs}

\usepackage{pdfcolmk}

\usepackage{graphicx}
\graphicspath{ {images/} }

\newcommand*{\rom}[1]{\expandafter\@slowromancap\romannumeral #1@}

\usepackage{xcolor}

\usepackage{wrapfig}
\newcommand{\tightunderset}[2]{%
  \mathop{#2}\limits_{\vbox to .3ex{\kern-0.95ex\hbox{$#1$}\vss}}}

\usepackage[utf8]{inputenc}

\usepackage{CJKutf8}

\usepackage{stmaryrd}
\usepackage{url}
\makeatother

\title {Mean Curvature in the Light of  Scalar Curvature.}

\author{Misha Gromov}

\begin{document}
\maketitle




We  think of   {\large \sf scalar curvature} as a Riemannin incarnation of {\large \sf  mean curvature}\footnote{Throughout the paper, we use the standard normalisation, where the unit sphere $S^{n-1}\subset \mathbb R^n$ has $mean.curv(S^n)=n-1$, the Ricci curvature of $S^n$ is $n-1$ and   the scalar curvature of $S^n$  is $n(n-1)$.}  and we search for  

 {\sf  constraints on  global geometric invariants of $n$-spaces $X$ with $Sc(X)>0$   that would generalise those for
 smooth {\it mean convex}  domains $X\subset \mathbb R^n$  also called domains with  
 {\it mean convex} boundaries $Y=\partial X$. i.e.  with $mean. curv(\partial Y)>0$.}\footnote{This is similar   in  spirit to parallelism   between spaces $X$  with positive sectional curvatures and convex subsets in $\mathbb R^n$, the best instance of which is Perelman's double sided bound on  the product of the $n$  {\it Uryson widths} of an $X$ by  $const^{\pm 1}_n\cdot vol(X)$.}\vspace {1mm}
 
And, as an unexpected bonus of this search, we   find out that \vspace {1mm}

{\sf  techniques developed for the study of manifolds $X$ with $Sc(X)\geq \sigma>0$ yield new
results for  
hypersurfaces $Y\subset \mathbb R^n$ with $mean.curv\geq \mu>0$.}

\vspace {1mm}

In what follows we briefly overview of what is known and what is unknown in this regard.  
\footnote {Part of this article  is a  slightly edited extract from my unfinished  paper {\sl 101 Questions, Problems and Conjectures around Scalar Curvature} \url{https://www.ihes.fr/~gromov/wp-content/uploads/2018/08/positivecurvature.pdf}}

  \vspace {2mm}
  \tableofcontents
  
\section{Inscribed Balls and Distance Decreasing Maps.} 

Let us recall classical  comparison theorems between radii of  balls in manifolds $X$ 
 with lower bounds on their Ricci curvatures by $\varsigma$ and on the mean curvatures of their boundaries\ by $\mu$ and 
 \vspace {1mm}
 
  \hspace{-6mm} {\sf the radii $R=R(n,\kappa,\mu)$ of the 
 the balls $B=B^n(\kappa,\mu)$ with $mean.curv(\partial B)=\mu$   in the standard (complete simply connected) $n$-spaces $X^n_\kappa$
 with sectional curvatures $\kappa=\varsigma/(n-1)$}, \footnote{If $\kappa\leq0$,  then the standard balls may only have  $\mu\geq (n-1)\sqrt {-\kappa}$.} \vspace {1mm}

  \hspace{-6mm} which go  back  to  the work by  Paul Levy, S. B. Myers and Richard L Bishop.  
\vspace {2mm}

{\sf Let $X$ be a (metrically) complete Riemannian $n$-manifold with a  boundary,  denoted $Y=\partial X$,  such that  $Ricci(X)\geq  \varsigma=(n-1)\kappa$, i.e. $Ricci(X)\geq  \varsigma\cdot g_X$,  
and  $mean.curv_\partial (Y)\geq \mu$.}

 For instance, $X$ may be a smooth, possibly unbounded  domain in  the Euclidean space $\mathbb R^n$, that is a closed subset bounded by a smooth hypersurface $Y\subset \mathbb R^n.$

\vspace{1mm}
  
   \textbf {Inball$^\mathbf n$-Inequality.}  {\sl The in-radius of $X$,    that is 
   $$Rad_{in}(X)=\sup_{x\in X} dist(x, \partial X),$$
    is bounded by the radius  $R=R(n,  \kappa, \mu)$  of the ball $B^n(\kappa, \mu)$ with
     $$mean.curv(\partial B^n( \kappa,\mu))=\mu$$ 
    in the standard $n$-space $X^n_\kappa$ of constant sectional  curvature $\kappa=  \varsigma/(n-1)$.}\footnote{If no such ball exists, e.g. if $\kappa\geq-1$ and $\mu=-(n-1)$ then, by definition, $R(n,  \varsigma, \mu)=\infty$.}
   
\vspace{1mm}

Indeed, the normal exponential map to  $\partial Y$ necessarily develops conjugate points on geodesic segments  normal to $\partial Y$ of length $>  R(n, \kappa, \mu)$.\vspace {1mm}
\vspace{1mm}

\textbf {Inball$^\mathbf n$-Equality} {\sl If $Rad_{in}(X)=R= R(n, \kappa, \mu)$, then $X$, assuming it is connected, is isometric to  an  $R$-ball in  $X^n_\kappa$.} \vspace {1mm}

This is proven by    fiddling at  the boundary points of the regions in $\partial X$,  where the maximal  in-ball meets $\partial X$.

\vspace{1mm}

\textbf  {$S^{\mathbf {n-1}}$-Extremality/Rigidity Corollary.}    \vspace{1mm}{\sl  Let $X\subset  X^n_\kappa$ be  a compact connected domain with smooth  connected boundary $Y=\partial X$.  }

{\sl Let  the mean curvature of $Y$ be  bounded from below by $\mu\geq 0$ and 
and let 
$$f:Y\to S^n=S^{n-1}(\kappa, \mu)=\partial B^n(\kappa, \mu)$$
be  a  distance non-increasing map  for  the distances in $Y  \subset  X^n_\kappa$ and in   $ S^{n-1}(\kappa, \mu)\subset X^n_\kappa$ induced from  that in the ambient (standard) space $X^n_\kappa$.}

  \vspace{1mm}

 \textbf {Ex$_{\sf\mathbf{mn}}$: Extremality of $S^{n-1}$}. {\sl If the map  $f$ is strictly  distance decreasing then $f$ is contractible.}\vspace{1mm}

\textbf {Ri$_{\sf\mathbf{mn}}$:  Rigidity  of $S^{n-1}$}. {\sl If $f$  is non-contractible, then 
it is an isometry.}\vspace{1mm}

{\it Proof.} Use  {\it Kirszbraun's theorem.}\footnote  {\url{https://en.wikipedia.org/wiki/Kirszbraun_theorem}} and extend map $f$ to a  distance non-increasing map  $F:X\to B= B^n(\kappa, \mu)$.

If $f$ has non-zero degree then so does $F$, hence,  the center of $B$ is in the image and the pullback of this center lies within distance $\geq  R(n, \kappa, \mu) $ from the boundary of $X$ and the above  {Inball$^\mathbf n$-inequality} and  {Inball$^\mathbf n$-equality} apply.

\vspace{1mm}

{\it Remarks} (a)  If $X$ is  isometrically realised by a domain   in the Euclidean or in the hyperbolic $n$-space and $dist(x_0, \partial X)\geq R$, then the normal  projection to the  ball,
 $X\to B_{x_0}(R)\subset X$, is a distance non-increasing map of degree 1.
 
 (b)  An essential drawback of  \textbf {Ex$_{\sf\mathbf{mn}}$}  and \textbf {Ri$_{\sf\mathbf{mn}}$} is an appeal to the {\it extrinsic metric}, that is the Euclidean  distance function restricted from   $\mathbb R^n$ to $Y$ and to  $S^{n-1}$, rather than to the {\it intrinsic metrics} in $Y$ and $S^{n-1}$ associated   with the {\it Riemannian metrics/tensors} induced from $\mathbb R^n$, where, observe,  
 the intrinsic metric in $Y$ may be incomparably greater than the extrinsic one.
 
Yet, we shall it in the  next section, 

\hspace {20mm}   {\it  \textbf {Ex$_{\sf\mathbf{mn}}$} remains valid   for the intrinsic metrics}

\hspace {-6mm}  but  the proof of this  relies 
on  {\it Dirac operators} on  manifolds  with {\it positive  scalar curvatures}    with no 
  no direct approach  in sight. (Intrinsic \textbf {Ri$_{\sf\mathbf{mn}}$} remains  {\it \large conjectural} at this point.)

\vspace{2mm}

The following refinements(s) of the inball$^\mathbf n$-inequality/equality   in the standard spaces of constant curvature, e.g in $\mathbb R^n$, must be also well known.  I  apologies to the author(s) whose paper(s) I failed to locate on the web.\vspace {1mm}

  \textbf {Inball$^{\mathbf n-1}$-Inequality}. Let  $X$ be a smooth  domain with $mean.curv(\partial X)\geq \mu$  in the standard $n$-space $X^n_\kappa$ with the sectional curvature $\kappa$ and let  $X$ contain a flat $R$-ball $B$ of dimension $(n-1)$, that is an $R$-ball in a totally geodesic  hypersurface in  $X^n_\kappa$.
  
  {\sl If $R\geq R(\kappa, \mu)$ then, in fact,  $R=R(\kappa, \mu)$ and $X$ is equal to an $R$-ball in $X^n_\kappa$.}

      \vspace {1mm}

{\it Proof.} Let $B_{ {\tightunderset {\smile}\frown}}(r)\supset B$ be the lens-like  region between two spherical  caps of hight $r\leq R$ and with the boundaries  $\partial B$. If $R>1$, the mean curvatures of these caps are $<\mu$; hence they do not meet $\partial X$  which makes 
the $R$-ball $B^n(R)= B_{ {\tightunderset {\smile}\frown}}(R)$ contained in $Y$  and the above two 
"inball$^\mathbf n$" apply.

{\it \textbf{ Intersection  Remark.}}  The above argument  also shows that if $X\subset \mathbb R^n$ with 
$mean.curv (\partial X)\geq n-1$ contains a  flat   $(n-1)$-ball,
$B_x^{n-1}(r)\subset X$, of radius $r<1$, then the distance from the center of this ball to the boundary of $X$ is bounded from below by    $dist(x, \partial X)> \frac {1}{2r^2}$. In words: 

\vspace{1mm}

{\sf if $X$ is {\it $r$-thin}  at a point $x\in X$,  then the  intersections of $X$ with hyperplanes passing through $x$  are {\it at most 
$\sqrt{2r}$-thick  at this point.}}\vspace{1mm}

In fact, the same apply to intersections of  $X$  with arbitrary hypersurfaces $S\subset \mathbb R^n$ where relevant constants  now  depend on the bounds on the  principal curvatures of $S$.

An instance of  such an  $S$, which  we shall meet  later in section 3,  is that of the sphere $S=S_x^{n-1}(1+\delta)$, $x\in X$, where  
\vspace{1mm}

{\sl the existence of an  $(n-1)$-disc of radius $r>20 \sqrt \delta$ in $S\cap X$ for this $S$    and      small $\delta>0$  necessities the existence  of a ball 
$$B_{x'}^n(2\delta)\subset X\cap B_x^{n}(1+\delta)$$ 
such that  $$dist(x', x)= 1-\delta.$$}

\vspace{1mm}

\textbf {Inball$^{\mathbf n-1}$ Mapping Corollary.} Let $X\subset \mathbb R^n$ be a  bounded domain with smooth  boundary  $Y=\partial X$, such that $mean.curv(\partial X)\geq n-1$ and let $ Y_+\subset Y$ be the intersection of $Y$
with the half  subspace $\mathbb R^n_+=\{x_1,...,x_i,...,x_n\}_{x_1\geq 0}\subset \mathbb R^n$. Then \vspace {1mm}

{\sl $ Y_+$ admits no proper distance decreasing map with non-zero degree to the unit ball $B^{n-1}(1)\subset \mathbb R^{n-1}$, where "distance decreasing" refers to the Euclidean distance  on $Y_+\subset \mathbb R^n$ and 
"proper" signifies  that  $\partial  Y_+\to S^{n-2}=\partial B^{n-1}(1)$.}\vspace {1mm}

{\it \textbf { Maximal Principle  for Principal Curvatures.}}  The maximum principle argument, which was  used for the proof of the  strict inball inequalities (but {\it not} their extremal  equality cases!) and their corollaries   for hypersurfaces $Y$ in $\mathbb R^n$ with $mean.curv(Y)\geq n-1$,  trivially extends to hypersurfaces  $Y\subset \mathbb R^n$, where the {\it maxima of principal curvatures are} $\geq 1$ at all points $y\in Y$, and the same applies to hypersurfaces in the standard spaces $X^n_\kappa$ with constant sectional curvatures $\kappa$ for positive and negative  $\kappa$ as well.
 
What is more interesting is the following, probably known, simple   generalisation of the above   {Inball$^{\mathbf n-1}$-mapping corollary} derived with  this "principle".   \vspace {1mm}

{\it \textbf {Sharp bound on the filling radius  in codimension 2.}} {\sf Let $Y_0^{n-1}\subset \mathbb R^n$ be a  smooth cooriented hypersurface with connected boundary $Z^{n-2}=\partial Y_0$,   such that the maxima of the  principal curvatures of $Y_0$ satisfy
$$ maxcurv(Y_0,y)\geq 1 \mbox {  for all $y\in Y_0$}.$$ }

{\it Then $Z^{n-2}$ bounds a submanifold\footnote {A priori, this $Y_0$ may be singular, but it can be made smooth for $ codim(Y_0)=1$.}  $Y_1^{n-1}\subset \mathbb R^n$, i.e. $\partial Y_1^{n-1}=Z^{n-2}$, such that
$$dist_{\mathbb R^n}(y, Z^{n-2})\leq 1  \mbox { for all  } y\in Y^{n-1}_1.$$}

Consequently,

{\it  $ Z^{n-2}$  admits no  distance decreasing map  with non-zero degree to the unit sphere  $S^{n-2}(1)\subset \mathbb R^{n-2}$ with non-zero degree, where "distance decreasing" refers to the Euclidean distances  restricted to  $Z^{n-2}\subset \mathbb R^n$ and to $S^{n-2}\subset \mathbb R^{n-1}$. }

\vspace{1mm}

{\it Proof.}  If $Z^{n-2}$ doesn't bound in its  $\rho$-neighbourhood in $\mathbb R^n$ then, by the {\it Alexander duality} there exists a simple  curve $C\subset \mathbb R^n$, with both ends going to infinity in $\mathbb R^n$,  which is {\it non-trivially  linked with $Z$} and such that 
$$dist_X(Z, C)> \rho.$$

We claim that there is a point $c\in C$ such that the $\rho$-sphere with the center $c$, say
$S_c^{n-1}(\rho)$ is tangent to $Y_0$ at some  point $y_0\in Y_0^{n-1}$, where, this sphere is  locally positioned  "inside" $Y_0$ relative to the given coorientation of $Y_0$.

To see this let $ \tau: S^{n-2}(\rho)\times C\to \mathbb R^n$ be the  map (tautologically) defined via the identification of the  sphere
 $S^{n-2}(\rho)=S_0^{n-2}(\rho)\subset \mathbb R^{n-1}$ with all $S_c^{n-1}(\rho)$ by  parallel translations and let us assume this map is transversal to $Y_0$.
 
 Then the $\tau$-pullback of $Y_0$,
 $$\Sigma=\tau^{-1}(Y_0)\subset  S^{n-2}(\rho)\times C$$
is a smooth hypersurface in $S^{n-2}(\rho)\times C$, the projection of which to  $S^{n-2}(\rho)$ has non-zero degree, namely the degree equal the linking number between $Z^{n-2}$ and $C$. 

It follows, that there is a  connected component, say $Z_0\subset Z$ which separates the two ends of the cylinder   $S^{n-2}\times C$, where one of the ends is regarded as "internal" with respect to the coorientation of $Y_0$ and the other one as "external". 

Then we let $(s,c)_0$ be  the minimum point of the function $(s,c)\mapsto c\in C=\mathbb R$  restricted to the "external" connected component of the complement to $\Sigma_0$ in  the cylinder $S^{n-2}(\rho)\times C$ and observe that
$$\tau((s,c)_0)\in Y_0\cap S_c^{n-1}(\rho)$$
serves as  the desired "internal kissing point" of  the sphere $S^{n-2}_c (\rho)$  with $Y_0$.

Finally, by the maximal principal, all principal curvatures of $Y_0$ at this "kissing point" are bounded by 
$\frac {1}{\rho}$  and  {\sf the proof is concluded.}

\vspace{1mm}

{\it Remark.}  The above  automatically generalises to   the standard spaces with constant curvatures and to many spaces with variable curvatures.

Also there is 
the following  (rather trivial)  version of this for submanifolds $Z^{n-k}\subset \mathbb R^n$ of codimensions $k\geq 3$: \vspace{1mm}

{\sl if $Z^{n-k}$ doesn't bound in its $\rho$-neighbourhood, and if  $Y_0^{n-k+1}\subset   \mathbb R^n$ has  $\partial Y_0^{n-k+1}=Z^{n-k}$, then there exists a  normal vector  $\nu_0$ to $Y_0$ at some point, such that all principal curvatures in the direction of $\nu$ at this point are bounded by $\frac {1}{\rho}$. }




\section {Manifolds  with Lower Bounds on their Scalar curvatures and  on the  Mean Curvatures of their Boundaries.}   

An essential link between positive  mean and positive scalar curvatures is furnished  by an elementary   observation  \cite {withLawson1980}
 that the natural {\it continuous} Riemannin metric on the double $X\cup_YX$ of a domain 
$\subset \mathbb R^n$ with boundary $Y=\partial X$ with {\it positive mean curvature} can be approximated by  smooth metrics with {\it positive scalar curvatures}.

\vspace {1mm}
This leads 
  to the following "intrinsic" improvement of  the  above "extrinsic  extremality" of the Euclidean spheres.

\vspace {1mm}

{\large \textbf \smiley}$^{ \hspace{-0.4mm}\to \hspace{-0.2mm}S^{n-1}}_{\hspace {-0.7mm}mn}$ {\it \textbf{ Sharp Bound on $\mathbf {Rad_{S^{n-1}}(Y\subset \mathbb R^n)}$}}.   {\it Let   $Y\subset \mathbb R^n$ be a closed hypersurface  and let $f$ be  a Lipschitz map of $Y$  to the unit sphere, 
 $$f:Y\to S^{n-1}=S^{n-1}(1)\subset \mathbb R^n.$$ 
If  $mean.curv (Y)\geq n-1$ and if   $f$ strictly decreases  the lengths of the  curves in  $Y$, then $f$ is contractible.} \vspace {1mm}

In fact,  this   is a corollary to the following  result derived  in \cite {boundary}  from  Goette-Semmelmann's estimates  \cite {GS2000}  for  {\it twisted Dirac operators}. 

 \vspace {1mm}

 {\large \textbf \smiley}$^{ \hspace{-0.4mm}\to \hspace{-0.2mm}S^{n-1}}_{\hspace {-0.7mm}}${\tiny $_{ \hspace{-8.5mm}Sc\geq 0}$}  \hspace {0mm} { \it \textbf{ Sharp Bound on $\mathbf {Rad_{S^{n-1}}(Y=\partial X)_{Sc(X)\geq 0}}$}}. {\sl If a Riemannin manifold  $Y$  serves as a  boundary with mean curvature $\geq n-1$ of a  compact Riemannian {\it spin} $n$-manifold $X$ with {\it  non-negative scalar curvature}, then $Y$ admits no non-contractible  map
$f:Y\to   S^{n-1}$ which strictly decreases  the lengths of curves in  $Y$.}

However, the following remains problematic.\vspace{1mm}

{\it \textbf {Question 1.}} {\sf  Is there a direct elementary proof of{\large \textbf \smiley}$^{ \hspace{-0.4mm}\to \hspace{-0.2mm}S^{n-1}}_{\hspace {-0.7mm}mn}$?}\vspace{1mm}

{\it \textbf {Question 2.}} {\sf Does {\large \textbf \smiley}$^{ \hspace{-0.4mm}\to \hspace{-0.2mm}S^{n-1}}_{\hspace {-0.7mm}mn}$  generalise to hypersurfaces $Y$ in the standard spaces $X^n_{\kappa}$ with sectional curvatures $\kappa\neq 0$?}
\vspace{1mm}

{\it \textbf {Question 3.}} {\sf Does the above  $\mathbf {Ri_{mn}}$, that is the rigidity   of $S^{n-1}\subset \mathbb R^n$,    hold (similarly to the  extremality   $\mathbf {Ex_{mn}}$) for the {\it intrinsic metrics} in $Y$ and $S^{n-1}$?}

Namely,

{\sl    are   maps of smooth closed hypersurfaces $Y\subset \mathbb R^n$ with  $mean.curv(Y)\geq n-1$ to spheres,  which do not increase the length of curves,
$$f:Y\to S^{n-1}=S^{n-1}(1),$$
 either  contractible or 
isometric?}
(This  may be significantly  easier  for {\it smooth} maps $f$.)\vspace{1mm}


{\it \textbf {Question 5.}} {\sf Is there a {\it sharp version} of the  {\textbf{ Inball$^{\mathbf n-1}$}-mapping corollary}  to the intrinsic metric in $Y$?}\vspace{1mm}

Namely, let $Y\subset \mathbb R^n$ be a smooth  closed hypersurface with mean curvature$\geq n-1$ and let $Z\subset Y$ be a hypersurface, $dim(Z)= n-2$, which divides $Y$ in two halves, say $Y_-\subset Y $ and  $Y_+ \subset Y$.\vspace{1mm}

{\sf Does  this  $Z=\partial Y_-=\partial Y_+$ admit a non-contractible  distance decreasing map to $S^{n-2}=\partial S_-^{n-1}=\partial S_+^{n-1}$, where the distance  in $Z$ is induced from the Riemannian distance in $Y$ and the distance in $S^{n-2}$ is the intrinsic spherical one (which is equal to the distance coming from  the ambient  sphere 
 $S^{n-1} \supset S^{n-2}$)?}\vspace{1mm}

Notice in this regard that  
{\large \textbf \smiley}$^{ \hspace{-0.4mm}\to \hspace{-0.2mm}S^{n-1}}_{\hspace {-0.7mm}mn}$ implies a non-sharp version of this, due to the following simple  corollary  of the  extension property 
of Lipschitz maps to $\mathbb R^n$.
\vspace {1mm}

{\sf \textbf{[Lip$_ {\hspace {-0.5mm}\mathbf{ \sqrt {n}}}$]}} \hspace {1mm}  {\sf  {\it 1-Lipschitz},  i.e.  distance non-increasing,  maps $Z=\partial Y_{\mp}\to \partial S_\mp^{n-1}\subset S^{n-1}$ of non-zero degree  (obviously) extend to 
$\lambda$-Lipschitz maps $Y\to S^{n-1}$ for $\lambda= \frac {\pi}{2}\sqrt{n-1}$ which also have non-zero degrees}.\vspace {1mm}

In fact, this extension property together with    {\large \textbf \smiley}$^{ \hspace{-0.4mm}\to \hspace{-0.2mm}S^{n-1}}_{\hspace {-0.7mm}}${\tiny $_{ \hspace{-8.5mm}Sc\geq 0}$}\hspace {2mm}  yield the following (alas, non-sharp) bound on the {\it (hyper)spherical radius} $Rad_{S^{n-2}}(Z\subset Y)$. \vspace {1mm}

{\large \textbf \frownie}$^{ \hspace{-0.5mm}\to \hspace{-0.2mm}S^{n-2}}${\tiny $_{ \hspace{-8.5mm}Sc\geq 0}$} \hspace {2mm}  {\sf Let  $X$ be a  compact $n$-dimensional  Riemannian spin manifold with boundary $Y=\partial X$, such that  $Sc(X)\geq 0$ and  $mean.curv(Y)\geq n-1$.}

{\sl If a  closed  hypersurface $Z\subset Y$ homologous to zero in $Y$ 
admits an 
$\varepsilon$-Lipschitz map $Z\to S^{n-2}=\partial S^{n-1}_\pm$ of non-zero degree, 
where "Lipschitz" is understood for the distance associated with the Riemannian metric in $Y$ induced from $X$, then 
$$\varepsilon \geq    \frac {2}{\pi \sqrt{n-1}}.$$}

\vspace{1mm}
{\it Remark/Example.} The "homologous to zero" condition is essential: non-contractible curves in  2-tori  $Y\subset \mathbb R^3$ with $mean.curv (Y)\geq 2$  may be uncontrollably long.

\vspace{2mm}

 Now  let us show that  a combination of the  above  argument  with that used in the proof of  the sharp bound for $fill.rad(Y)$ in section 1 yields the following.

\vspace{1mm}

{\it \textbf  { Filling Radius Bound for $Z^{n-2}\subset \mathbb R^n$}.} {\sf Let $Y\subset \mathbb R^n$ be a smooth closed  (embedded without self-intersection!)  connected hypersurface with $mean.curv(Y)\geq n-1$ and let $Z=Z^{n-2}\subset X$ be a closed  (embedded without self-intersection) 
hypersurface in $Y$}.

{\sl If $Z$ is homologous to zero in $Y$ then it is homologous to zero in its $\rho$-neighbourhood in $\mathbb R^n\supset Z$ for
$$\rho\leq  \frac {\pi }{2}\frac {n}{n-1}.$$}
 

 {\it Proof.} If $Z$ doesn't bound in its  $\rho$-neighbourhood in $\mathbb R^n$ then, by the {\it Alexander duality} there exists a  closed curve $S\subset \mathbb R^n$, which is {\it non-trivially  linked with $Z$} and such that 
$$dist_X(Z, S)> \rho.$$
Then the radial projections of $Z$  to the $\rho$-spheres with the centers $s\in S$, say
$$f_s:Z\to S_s^{n-1}(\rho),$$
are {\it distance decreasing},  for the {\it Euclidean} distances restricted to  $Z$ and to $S_s^{n-1}$  
while  the resulting  map from $Z\times S $ to the $\rho$-sphere, call it
  $$f:Z\times S \to S^{n-1}(\rho)= S_{\mathbf 0}^{n-1}(\rho),$$
has non-zero degree, since this degree is equal to the linking number between $Z$ and $S$.

Then, using Kirszbraun theorem and an obvious $\frac {\pi}{2}$-Lipschitz homeomorphisms from the unit ball onto the unit hemisphere, $B^n\to S^n_\pm$, one extends this  $f$ to a $\frac {\pi}{2}$-Lipschitz map 
$$F: Y\times S\to S^n(\rho),$$
such that  
 $deg(F)=deg(f)\neq 0$, where, by using (the corresponding version of) {\sf \textbf{Lip$_ {\hspace {-0.5mm}\mathbf{ \sqrt {n}}}$}}
  and by  making  the curve $S$  longer if necessary,  one gets such an $F$ with
$$Lip(F)\leq  \lambda =\frac {\pi}{2}.$$

 Finally  observe that $Y\times S$ serves as the boundary of the  manifold $X^{n+1}$, that is the domain in $\mathbb R^n$ bounded by $Y$  times $S$, where $Sc(X)=0$ and $mean.curv(Y\times S)=mean.curv(Y)\geq n-1$  and that 
  {\large \textbf \smiley}$^{ \hspace{-0.4mm}\to \hspace{-0.2mm}S^{n}}_{\hspace {-0.7mm}}${\tiny $_{ \hspace{-6mm}Sc\geq 0}$} 
 \hspace {2mm} implies that 
 $$ \rho\leq \frac{\lambda n}{n-1},   $$
by which {\sf the proof is concluded. }

\vspace {1mm}

{\it Remarks and Conjectures.} (a) Our inequality with   $const= \frac {\pi }{2}\frac {n}{n-1}>1$
is not especially exciting in view of  the sharp  estimate (where $const=1$) available  with  the    maximum principle. 
But the above argument applies to   a class of  spaces quite different from $\mathbb R^n$, where   the scalar curvatures are bounded from below, and which admit (exactly or approximately) contracting projections to the balls, such, for instance as 
  {\sf manifolds  with $Sc\geq -\varepsilon $ which are  bi-Lipschitz equivalent to $\mathbb R^n$}. \vspace {1mm}

However,  {\it sharp} filling inequalities  in this kind of  spaces remain conjectural.\vspace {1mm} 

 (b) Let  a hypersurface $Z\subset  Y=\partial X$, where $X$ is a compact Riemannin manifold, divide $Y$ into two halves, say  $Y_+$ and $Y_-$, such that 
 $$mean.curv(Y_+)\geq n-1\mbox {  and  } mean.curv(Y_1)\geq -(n-1). \leqno {[\mathbf  {MN_{\pm(n-1}}]}$$
 Then one knows  that $Z$ bounds a  {\it stable hypersurface}  $Y_{min}\subset X$   with {\it constant mean curvature} $n-1$  (called "$\mu$-bubble"  
 in  \cite {positive}).\footnote{If $If n\geq 9$, this $Y_{min}$ may, a priori, have "stable singularities"  but these  can be rendered harmless  with the recent Schoen-Yau theorem in \cite {SY2017} and/or, possibly, using  the wotk  by Lohkamp \cite {L2018}.} 
 \vspace {1mm}
 
If $Sc(X)\geq 0$,  this  yields -- we show it somewhere else --   alternative proofs of

{\sf non-sharp bounds on the (hyper)spherical radius and   on the  filling radius of  $Z$  provided  the mean curvature of $Y\supset Z$ satisfy the condition $[\mathbf  {MN_{\pm(n-1}}]$}, which is weaker than $mean.curv(Y)\geq n-1$.\vspace {1mm}



 \vspace {1mm}

(c)   
 Let $X^+$ be a complete Riemannian $n$-manifold, which, for simplicity's sake, we assume having uniformly bounded local geometry, such as $X^+=\mathbb R^n$, for instance.

Let $Z=Z^{n-2} \subset  X^+$ be a smooth closed oriented submanifold which  bounds a  cooriented submanifold $Y_0=Y^{n-1}_0\subset X^+$ with $mean.curv (Y_0)\geq n-1$.

   \vspace {1mm}

{\sf Does then $Z$ bound a stable  hypersurface $Y_{min}\subset X^+$  with $mean.curv(Y_{min})=n-1$?}   \vspace {1mm}

( Intuitively, this $Y_{min}$  could be  obtained by an "inward deformation"  of $Y_0$, where, in  general, even if $Y_0$ is compact, the bubble $Y_{min}$ may somewhere go to infinity; yet, this is  ruled out by  the "bounded geometry" condition    as in  11.6 in \cite {inequalities}).

If  this works, we would  obtain   bounds on  filling (hyperspherical)  radii and filling volumes of $Z$, in manifolds $X^+$  with $Sc(X^+)\geq \sigma$,
 similar to the   bound in $\mathbb R^n$ obtained  with the maximum principle.

But since  minimal bubbles $Y_{min}$  can  intersect  $Y_0$  away from $\partial Y_0 =\partial Y_{min}$,  the existence  of  these $V_{min} $ seems  problematic. 
\vspace{1mm}

(d) Let  $Z=Z^{n-2}\subset X^+$, e.g. for $X^+=\mathbb R^n$ be as above.\vspace{1mm}

Can  our bound on  $fill.rad(Z)$ be upgraded  to such a  bound for   an   $(n-1)$-volume minimizing filling of $Z$ that is a  volume minimizing  hypersurface, say   $Y_\ast \subset X^+$ with $\partial Y_\ast=Z$?

\vspace{1mm}

Namely 

{\sf Is there a bound on 
$\sup_{y\in Y_\ast}  dist (y, Z)$ for our $Z$}? in terms of the lower bounds on the scalar curvature of $X^+$
and on the mean curvature of some $Y\supset Z$?
\vspace{1mm}

(Ideally we would like this "{\it dist}"  to be associated with the Riemannin metric in $Y_\ast$ rather 
than with that in $X^+$.)

\vspace{1mm}




\section {Back to Mean Convex Domains  in $\mathbb R^n$}

The following may   be known but I couldn't find a reference.   \vspace {1mm}

 {\it \textbf {Conjecture 1. }} {\sf Let  $X$ be an (infinite!) {\it mean convex} domain  $X\subset \mathbb R^n$, i.e. $mean.curv(\partial X)\geq 0$. 
 
 If the boundary  $\partial X$ is disconnected, then none of the connected component of $\partial X$  may have  its mean curvature separated away from $0$, i.e. infima of the mean curvatures of all components are zero.}
 
  \vspace {1mm}

This is known in the case   $n=3$, where the following   better result is available.

 \vspace {1mm}

 \textbf{ Strong Half-space/Slab 
Theorem} \cite{Lopez}.
 {\sf The only mean convex domains,  i.e. with $mean.curv\geq 0$,   in $\mathbb R^3$ with disconnected boundaries are slabs between parallel planes.}  \vspace {1mm}
 
In fact, this easily follows  from \vspace {1mm}

 \textbf { Fischer Colbrie-Schoen Planarity   Theorem.}
{\sl Complete stable minimal surfaces in $\mathbb R^3$ are flat.} \vspace {1mm}

{\large \it Warning.}  The  Euclidean  spaces 
$\mathbb R^n$ for $n\geq 4$, contain  mean convex  domains bounded by pairs of   $n$-dimensional catenoids   \cite {Hoffman}
\vspace {1mm}

Mean convex  domains $X\subset \mathbb R^3$, even  with connected boundaries,  are subjects to strong geometric constrains.\vspace {1mm}

 {\large \LEFTCIRCLE}\hspace {-1.1mm}-Example.   {\it If  a  mean convex $X\subset \mathbb R^3$     contains a plane, then, assuming $\partial X$ is  non-empty and  connected,  $Y$  is equal  to a half-space.} \vspace {1mm}
 
   This follows from    {\it the  half-catenoid maximal  principle}  that was originally used by  Hoffman  and Meeks  \cite {Hoffman}  to  show that 
    \vspace {1mm}

\hspace {14mm}  {\it properly embedded   minimal surfaces $Y\subset \mathbb R_+^3$  are flat.}
    \vspace {1mm}

It may be  unclear what to expect  of general mean convex domains $X\subset \mathbb R^n$  but if 
the mean curvature  is separated away from zero, then one expects the following.    \vspace {1mm}

{\it   \textbf {Conjecture 2.}} {\sf If a domain $X\subset \mathbb R^n$ has
 $mean.curv(\partial X)\geq n-1$, then there exists  a continuous self-map
$R: X\to X$, such that

$\bullet$ \hspace {1mm} the image $R(X)\subset X$ has topological dimension 
$n-2$;}\vspace {1mm}

$\bullet$ \hspace {1mm} $dist(x, R(x))\leq const_n$ for all $x\in X$, 
with the best expected $const_n=1$.\vspace {1mm}

 {\it Comments.} This is 
 a baby version of the corresponding {\it \textbf {conjectural}} bound on the {\it Uryson width}  of complete (also non-complete?) $n$-manifolds $X$ with $Sc(X)\geq n(n-1)$: 
 
 {\sf there exists a continuous map from $X$ to an $(n-2)$ dimensional, say   polyhedral, space, say
 $\rho:X\to P^{n-2}$, such that the diameters (also areas?)  of  the pullbacks of all points $p\in P$ are bounded by
 $$diam( \rho^{-1}(p))\leq  const_n,$$
 possibly  with $const_n=1$ (and $const_n=4\pi$ in the case of areas).}\footnote {See \cite {withLawson} and  
 \cite {MN 2011]} for the case  $n=3$.}

 Ideally, one would like to have the above  map  $R$ as the end   result of  a  homotopy, a kind of mean curvature flow, which would collapse $X$ to $P$ and blow up the mean curvature at all points in the process.
 
 Similarly,  $\rho$ may result  from a Ricci kind of flow
 which would shrink $X$ to $P^{n-2}$ in finite time with a simultaneous blow up of the scalar curvature.

\vspace{1mm}\vspace {1mm}

{\it Exercises.} (a) Let $n=3$  and show that there exists a map $f:\partial X \to \mathbb R^3$, such that  the image $R(X)\subset X$ has topological dimension 
$1$  and  $dist(x, R(x))\leq 100$ for all $x\in \partial X$. 

{\it Hint.} Argue as in the proof of corollary  10.11 in \cite {withLawson}.\vspace{1mm}

(b) Let $X$ be a complete $n$-manifold with {\it disconnected} boundary, where
 $mean.curv(\partial X)\geq \mu\geq 0$. If $n\leq 7$. Then, for all positive $\nu\leq \mu$, $X$ contains an  $n$-submanifold  $X_\nu $  with a  {\it disconnected} boundary  which have  {\it constant} mean curvature $mean.curv( \partial X_\nu)=\nu$. 
\vspace{1mm}

{\it Remark.} This may be only of "negative" use for  $X\subset \mathbb R^n$, where it may help to  
settle  conjecture  2, but it may be more relevant for domains in the hyperbolic space $\mathbf H^n$ with the  sectional
curvature $-1$,  where mean convex domains of all kinds are abundant  and  where the counterpart of conjecture 2 refers to $X\subset \mathbf H^n$  with $mean.curv ( \partial X)>n-1$.  (Complete non-compact hypersurfaces with 
$mean.curv =n-1$ in $\mathbf H^n$, which have no Euclidean analogues,  look especially intriguing.)

Also it may be amusing to look at the conjecture 2 from the position of manifolds with $Ricci\geq 0$, where the natural  guess  is  as follows.\footnote{Another natural guess is that the  answer to this,  along with conjecture 4, must be known to right people.}\vspace {1mm}

{\it   \textbf {Conjecture 3.}}  {\sf If an $X$ with $Ricci (X)\geq 0$ has   disconnected boundary with the mean curvature $\geq -\mu$,  then no connected component of this boundary  can have mean curvature $\geq \mu+\varepsilon.$}

(This is, of  course, obvious for compact $X$.)

\vspace {1mm}\vspace {1mm}

\hspace {43mm}{\large \it\textbf {Symmetrization.}}  

\vspace {1mm}

Intersections of mean convex subsets  $X$ in Riemannin manifolds are mean convex   with a properly defined  generalised mean curvature and, more generally, the inequality $mean.curv( \partial X)\geq \mu$, is stable under finite and infinite  intersections of such domains for all $-\infty < \mu <\infty$. 

 This allows {\it $G$-symmetrization} of $X$'s under actions of isometry groups $G$ acting on the ambient manifold $X^\ast\supset X$, e.g. for $X^\ast=\mathbb R^n\supset X$
 $$X\leadsto X_{sym}=\bigcap_{g\in G}g(X).$$

This suggests, for instance,   the  proof of the above   {\large \LEFTCIRCLE}   in some (all?) cases  by symmetrization  $X\leadsto X_{sym}\subset \mathbb R^3$  where  $X_{sym}$ is equal to the   intersection of the copies of $X$ obtained by rotations of  $X$ around the axis normal to the hyperplane   $\mathbb R^{2} \subset X$. \vspace {1mm}
 
Similarly, one  can apply symmetrizations to   mean convex domains $X\subset \mathbb R^n$ for $n\geq 4$ which contain hyperplanes  and  where definite results need fast decay conditions on the mean curvatures 
of $\partial X$.
For instance, one shows with symmetrization -- this must be classically  known --  that \vspace {1mm}

{\sf The hyperplane $\mathbb R^{n-1}\subset \mathbb R^n$ admits no mean convex perturbations with {\it compact supports,}}

 \vspace {1mm}

  \hspace {-6mm}which also can be derived,  from non-existence of  
$\mathbb Z^n$-invariant metrics with $Sc>0$ on $\mathbb R^n$.

In fact, such  perturbations admit  $\mathbb Z^{n-1}$-invariant extensions, more precisely invariant under the action of the group (isomorphic to  $\mathbb Z^{n-1} $) generated by sufficiently large mutually normal translations. and at the same time the metric in $X$ can be perturbed to have the scalar curvature $>0$. Then the metric on the double $X\cup_\partial X$ of $X$ can be similarly made  $\mathbb Z^n$-invariant keeping $Sc>0$.

  \vspace {1mm}

We mentioned all this   to  motivate, albeit not very convincingly,  the following.  \vspace {1mm}

 {\it  \textbf {Conjecture 4.}}  {\sf The only $\mathbb Z^{n-3}$-invariant  mean convex domains in $\mathbb R^n$ with disconnected boundaries are slabs between parallel hyperplanes.}  \vspace {1mm}

  \vspace {1mm}

 {\it  \textbf {Vague Question.}} {\sf What could be  a counterpart of the mean curvature symmetrization for metrics with $Sc\geq \sigma$?}
  \vspace {2mm} 

{\it \textbf {Stability of the Inball$^n$-Inequality.}} Let the boundary of a smooth domain $X\subset \mathbb R^n$
have
$mean.curv (\partial X)\geq n-1-\varepsilon$ for $\varepsilon >0$ and let $x_1,x_2\in X $  satisfy
$dist(x_1,\partial X)\geq  1$, and $dist(x_2,\partial X)\geq r$ for $0<r\leq 1$

Then 
$$\mbox {either   $  dist(x_1,x_2)\leq 1-r+\delta$ or  $dist(x_1,x_2)\geq 1+r-\delta$},$$
where $\delta\leq \delta(\varepsilon)\to 0$ for $\varepsilon \to 0$.

  \vspace {1mm} 
  
  {\it Proof.} Symmetrize around the axes between $x_1$ and $x_2$.

  \vspace {1mm} 

{\it Corollary.}  If $\varepsilon$ is small, then  the {\it unit balls} which are  contained in $X$ {\it can't be continuously moved}  $\delta$-far in $X$ from their original positions. 

Moreover, the part of $Y=\partial X$ which lies $\delta$ close to the ball $B^n_{x_1}(1)\subset X$ has only small holes, about $\sqrt\delta $ in size, " which can be sealed by a small perturbation of $Y$ and thus lock  $B^n_{x_1}(1)$ in the concentric ball of radius $1+\delta$. This means, more precisely, the following.

\vspace {1mm}

 {\sf There exists    continuous self-mappings  $f= f_\varepsilon:X\to X$, $\varepsilon>0$, such that}  \vspace {1mm}

$\bullet_1$   {\sl $f$ is  supported (i.e. $\neq Id$) in a given arbiltrarilly small neighbourhood of the boundary $Y=\partial X$,

$\bullet_2$   $f$  moves all points  by  a  small amount for small $\varepsilon$,
 $$dist(f(x), x)\leq  \delta(\varepsilon)\underset{\varepsilon\to 0}\longrightarrow 0, $$

  $\bullet_3$  the image $f(Y)\subset X$ "locks"  all unit balls in $X$, that is 

{\sf every point $x\in X$ with $dist(x,\partial X)\geq 1$ is contained in a connected component of the complement $\mathbb R^n\setminus f(Y)$, where this components itself is contained in the ball $B_x(1+\delta)$, such that  this $\delta$ 
is also bounded by   $\delta(\varepsilon)\underset{\varepsilon\to 0}\longrightarrow 0.$}} 

\vspace {1mm}

{\it Proof.} Confront the above with the intersection corollary from section 1 and conclude that if  the unit ball $B_x(1)$
is contained in $X$,  then the intersection of $X$ with   the  $\delta$-greater   concentric sphere $S=S^{n-1}(1+\delta)\subset \mathbb R^n$ contains no disc of radius $>20\sqrt{\delta}$

Then collapse this intersection, call it  
$$V=X\cap S\subset S=S^{n-1}(1+\delta),$$
 to its cut locus  with respect to the boundary, denoted $\Sigma\subset  V$, and then extend the map $V\to \Sigma$
to the required map $X\to X$.\vspace{1mm}

{\it \textbf {Stability in the Limit.}} Given a domain $X\subset \mathbb R^n$  let $\Delta_X$ be the distance function to the boundary inside $X$ extended by zero outside $X$,
$$\mbox {$\Delta_X(x)= dist(x,\partial X)$ for $x\in $X and $\Delta_X(x)=0$   for $x\in \mathbb R^n\setminus X$.}$$

Say that a sequence of  subsets  $X_i\subset \mathbb R^n $ {\it regularly converges} if the functions $\Delta_{X_i}$ uniformly converge on compact subsets in $\mathbb R^n$.  Then the {\it regular limit} $X_{\infty^\circ}$ of $X_i$ is defined  as the (open!) subset, where the limit function  $\lim_{i \to \infty}\Delta_{X_i}$ is {\it strictly positive.}

(Observe  that every sequence $X_i$ contains a regularly convergent subsequence.) 

Now, the existence of the above  $f:X\to X$ (obviously) implies  the following.\vspace {1mm}

{\sl If $mean.curv(\partial X_i)\geq \mu_i \underset {i\to \infty}\longrightarrow n-1$, and $X_i$ regularly converge to $X_{\infty^\circ}$,
then the connected components in $X_{\infty^\circ}$ of the points $x\in X_{\infty^\circ}$ where $dist(x,\partial X_{\infty^\circ} )=1$ are exactly  open unit balls with the centers at these points $x$.}\vspace {1mm}

{\it Remark.} The non-inclusion $ dist(x_1,x_2) \notin [1-r+ \delta, 1+r-\delta]$ is sharp:   \vspace {1mm} 

{\sl if $dist (x_1,x_2)>2$, then, for all $\varepsilon>0 $, there are  smooth domains $X$ in $\mathbb R^n$ with $mean.curv(\partial X)\geq n-1-\varepsilon$, which contains both  unit balls $B_{x_1}(1)$ and  $B_{x_1}(1)$,}
as we shall explain in section 5.  \vspace {1mm}


{\it \textbf {Conjecture 5.}}\footnote{This is motivated by an aspect of the Penrose conjecture explained to me by  Christina Sormani.} Let $X_i\subset \mathbb R^n$  be a sequence of smooth  domains, such that all of them contain the unit ball $B=B^n_0(1)$ and such that
$$mean.curv(X_i)\geq \mu_i\to n-1 \mbox {  for } i\to \infty.$$
Then 

{\sf there exists a sequence of compact  domains $B'_i$ with smooth boundaries $S'_i=\partial B'_i $, which approximate $B$ from above, i.e. 
 $$\bigcap_iB'_i= B,$$
and such that 
$$vol_{n-1}(S'_i\cap X_i)\to 0\mbox {  for }i\to \infty.$$}

{\it Exercises.} (a) Decide what are   possibilities for the values of the  distance between $x_1$ and $x_2$, in the case 
$mean.curv(\partial X)\geq 1$ and 
$dist(x_i, \partial X)=\frac{n-2}{n-1}+\varepsilon_i$, $i=1, 2$.\vspace {1mm}

(c)\footnote {This is a replacement to  the {\it  erroneous}  
\PlusThinCenterOpen-inequality  from my "101-Questions" paper.} Let 
$$X_0=B^{n_1}(r_1)\times B^{n_2}(r_2) \subset \mathbb R^{n_1+n_2}$$
and evaluate the maximal $\mu$, for which there exists  $X\supset X_0$  with $mean.curv(X)\geq \mu$.

 \vspace {1mm}
{\it Hint/Remark.} $G$-Symmetrization where $G$ is product of two  orthogonal groups, $G=O(n_1)\times  O(n_2)$,  renders the problem 1-dimensional,  the analysis of which -- I haven't tried it myself -- seems easy. 

But finding this maximal $\mu$  for products of  $k$ balls, 
$$X_0=B^{n_1}(r_1)\times ...\times B^{n_k}(r_n)\subset \mathbb R^{n_1+...+n_k},$$ where $k\geq 3 $, may be more difficult.  
\vspace{1mm}

(d)  {\it $(\delta, \mu)$-Regularisation.} 
 Given a  closed (i.e the closure of open)  domain $X$ in a  complete  Riemannin manifold $X^\ast$, 
 let 
$X_{\pm_\delta}\subset X$  be the $\delta$-neighbourhood of the 
$\delta$-sub-level   $\Delta_X^{-1}[,\delta]\subset X$, that  
 is the union of the $\delta$-balls from  $X^\ast$ which are contains in $X$.

Then let $X_{{\pm_\delta},\mu}$ be the intersection of all subsets in $X^\ast$ which contain $X_{\pm_\delta}$
and have the mean curvatures of their boundaries $\geq \mu$. Show that:

$\bullet_1$ the operation $X\mapsto X_{\mp_\delta}$ is {\it idempotent}, moreover,
$$ (X_{\mp_\delta})_{\mp_{\delta'}}=X_{\mp_\delta}\mbox {   for  }\delta'\leq \delta;$$

$\bullet_2$ if $X$ is compact  mean convex with a  {\it piecewise smooth} boundary, then the boundary $\partial X_{\mp_\delta}$ is  {\it $C^1$-smooth} for
small $\delta>0$ and 
$$X_{{\pm_\delta},\mu}=X_{{\pm_\delta}}.$$
for  $ \mu<<\frac {1}{\delta}.$\vspace {1mm} 

{\it Remark.} The operation $X\mapsto X_{\mp_\delta}$ with (relatively) large $\delta$ doesn't seem to preserve   mean convexity even in $\mathbb R^n$, for $n\geq 3$, where the apparent example can be obtained -- I didn't check it all 100\% -- by a $C^\infty$-small perturbation of  $X_0=\blacktriangle\times \mathbb R^{n-2} \subset \mathbb R^n$, where one of the angles of the triangle $\blacktriangle\subset \mathbb R^2$ is very small.
 
 This raises the problem of evaluating the maximal distance from $X_{{\pm_\delta},\mu}$ to $X_{\pm_\delta}\subset X_{{\pm_\delta},\mu}$ for $\mu$-convex domains $X$,  where this "maximal  distance" is
$$\sup_{x \in X_{\pm_{\delta,\mu}}}dist (x,X_{\pm_\delta}).$$.


\section {Being Thick in Many Directions}

Start by observing that the products of  balls $B^{n_1}=B^{n_1}(1)\subset \mathbb R^{n_1}$ by Euclidean spaces,
 $$B^{n_1}\times  \mathbb R^{n_2}$$ 
 are {\it mean curvature extremal domains} in $\mathbb R^{n_1+n_2}$:
 
 {\sl if a domain $X\subset \mathbb R^{n_1+n_2}$ 
 with the mean curvature $\partial X\geq n_1+n_2-1$ contains $B^{n_1}\times \mathbb  R^{n_2}$
 then 
 $$X=B^{n_1}\times \mathbb   R^{n_2}.$$}
 
This can be shown either by $G$-symmetrization of $X$ for $G$ being the isometry group of 
 $B^{n_1}\times  \mathbb R^{n_2}$ or by an application of the maximum principle to a family  of domains 
 $$X_t\supset X_0=B^{n_1}\times  \mathbb R^{n_2},$$
  for 
 $$X_t= \bigcup_{x\in  \mathbb R^{n_2}}  B^{n_1}(1+\phi_t(x))\subset \mathbb R^{n_1+n_2}=\mathbb R^{n_1}\times \mathbb R^{n_1+n_2},$$
where $$\phi_t:  \mathbb R^{n_1+n_2}\to\left [1, 1+\frac {n_1+n_2-1}{n_1-1}\right]$$
are smooth functions, such that 
$$\phi_0=1, \hspace {1mm}  \phi_1= \frac {n_1+n_2-1}{n_1-1} \mbox {   }\frac {d\phi_t}{dt}>0,$$ and where
 the first and the second partial derivatives of $\phi)t(x)$ are very small and fast 
 decreasing as $t$ grows.
 \vspace {1mm}

\vspace {1mm}

The extremality of $B^{n_1}\times  \mathbb R^{n_2} $ motivates the following conjectural bound  the {\it macroscopic dimensions} of the sets of (the centers of) large balls in $X$, where,
recall,   

 {\sf the  macroscopic dimension  of  a metric space  $M$ is the minimal dimension of 
a polyhedral space $P$, for which  $M$ admits a continuous map $\Delta:M\to P$, such that
$diam_M(\Delta^{-1}(p)) \leq d$  for some constant  $d=d(M)<\infty$.}

\vspace {1mm}

 {\it \textbf {Conjecture 6.}}  {\sf Let a domain  $X\subset \mathbb R^n$ satisfy 
 $$mean.curv(\partial X)\geq n_1-1+\varepsilon $$  
 for   $n_1<n$ and  $\varepsilon>0 $.  Then the macroscopic dimension of the   subset  $X_{-1}\subset X$ of the points $x\in X$, such that $dist(x,\partial X) \geq 1$, satisfies:
 $$macr.dim(X)< n-n_1,$$
 where the equality  $macr.dim(X)=n-n_1$ with $\varepsilon =0$ is achieved {\it only} for 
 $$X_0= B^{n_1}\times \mathbb R^{n-n_1}.$$}

{\it Remarks.} (a)   {\it The case of $n_1=n-1.$}  Conjecture 6 in this case says that if  $mean.curv(\partial X)\geq n-1-\alpha$
for $\alpha <1$, then the diameters of all connected components of   the above subset $X_{-1}\subset X$ are bounded by $\beta \leq \beta(\alpha)<\infty$. 

This can be shown  for  $\alpha\leq \alpha_n\sim \frac {1}{n}$   by the  argument used for proving  stability  of the Inball$^n$-inequality,  where  slightly  better estimates on $\alpha_n$ can be, probably, obtained  with symmetrization  around $k$-planes containing centers of suitable $(k+1)$-tuples of unit balls in $X$.


\vspace {1mm}

(b)  The macroscopic dimension of $B^{n_1}\times \mathbb R^{n_2}$ is equal to $n_2$   by \vspace {1mm}

{\it Lebesgue's Lemma}: {\sf A continuous map $[0,1]^n\to \mathbb R^n$ necessary brings a pair of points from {\it opposite faces} in the cube $[0,1]^n$ to {\it  a single  point} in $\mathbb R^n$}.
\vspace {2mm}


\vspace {1mm}\vspace {1mm}

There is a  counterpart to   conjecture $6$ in the context  of Riemannian 
  manifolds  (with and  without boundaries) which express the idea 
that  \vspace {1mm}

{\sl if the scalar curvature of an $n$-dimensional  manifold $X$ is bounded from below by $n(n-1)$ then  the space of "large"  balls $B_x(r)\subset X$, say  of radii 
$r\approx  \pi$,  must be small, where the size of a ball is evaluated in comparison with geodesic  balls in the unit sphere $S^n$.}\vspace {1mm}

A  conceptually simple  instance  of  this concerns maps of closed $n$-manifolds $X$ to the unit sphere $S^n$, namely
\vspace {1mm}

{\sf the  space   of distance decreasing maps of non-zero degrees  $X\to S^n$, which we denote 
$ Lip_1(X_{\neq 0}^{\to\circ})$,}
\vspace {1mm}

 \hspace {-6mm} and  which we
  endow  with \vspace {1mm}
  
  {\sf the metric associated in the usual way  to the natural {\it length structure} in the space of maps to $S^n$,}
   \vspace {1mm}

   where the length of a curve in this space, that is a family of maps $f_t:X\to S^n$, is defined as the supremum of the lengths of the $t$-curves in $S^n$ drawn by   individual points $x\in X$,
$$length(f_t)=\sup_{x\in X} length (f_t(x)).$$

{  \it \textbf {Conjecture 7.}} {\sf  If $Sc(X)\geq n_1(n_1-1)+\varepsilon$, $\varepsilon>0$, then  
 $$macr.dim(Lip_1(X_{\neq 0}^{\to\circ}))< n-n_1.$$}

The simplest instance of this conjecture concerns manifolds $X$ with  $Sc(X)\geq n(n-1)-\varepsilon$, 
where it says that all connected components of the space $Lip_1(X_{\neq 0}^{\to\circ})$ are bounded. In fact, when $\varepsilon \to 0$, one expects more of this space, which we formulate as follows.

\vspace {1mm}

{\it \textbf {Conjecture 8:}} {\it {\textbf{ Stability of $S^n$.}}  {\sf Let $X$ be a closed  orientable Riemannin $n$-manifold,  such that $Sc(X)\geq n(n-1)-\varepsilon$. Then the diameters of the connected components of  the  quotient space  of
$ Lip_1(X_{\neq 0}^{\to\circ})$ under the orthogonal group action on $S^n$ satisfy:
$$diam (conn (Lip_1(X_{\neq 0}^{\to\circ})/O(n)))\leq \delta\leq \delta(\varepsilon)\underset {\varepsilon \to 0}\longrightarrow 0.$$}}

{\it Digression.}  To get a feeling for our metric in $ Lip_1(X_{\neq 0}^{\to\circ})$ observe that  such length   metrics  are defined  for spaces  maps to  all  length metric spaces $S$  and look at a few examples.

\vspace {1mm} 
 
 (i) {\it Continuous Maps}.   If $S$ is a  compact locally contractible space   and $B^n$ is the topological $n$-ball, then the  space $C(B^n$\hspace {-0.5mm}$\to$$S)$ of continuous maps $B^n\to S$, $n\geq 1$, has finite diameter if and only if $S$ has finite fundamental group. 

For instance  
$$\mbox {$diam( C(B^n$\hspace {-0.5mm}$\to$$S^m(1)))$$ =diam (S^m(1))=\pi$  for $m>n$}. $$
Somewhat less obviously,
$$\mbox {$diam( C(B^n$\hspace {-0.5mm}$\to$$S^n(1)))$$\leq 3\pi$}, $$
which implies that the ({\it infinite} cyclic) universal covering of the space
$C(S^{m-1}$\hspace {-0.5mm}$\to$$S^m(1))$ also has diameter $\leq 3\pi$.

More generally, 
$$\mbox {$diam( C(B^n$\hspace {-0.5mm}$\to$$S))$$\leq n\cdot const(S)$  }, $$
for all compact, say cellular, spaces $S$ with finite fundamental groups  and,
probably,  the universal coverings of all connected components of the spaces  $C(X$\hspace {-0.5mm}$\to$$S)$ are similarly  bounded
by $dim(X)\cdot const(S)$.




The above  linear bound on $diam( C(B^n$\hspace {-0.5mm}$\to$$S))$
is asymptotically matched by a  lower bound  for most (all?) compact  non-contractible spaces $S$.

  For instance,  it follows from 1.4 in \cite {dilatation}  that 

$$\mbox {$diam(C(B^n$\hspace {-0.5mm}$\to$$S))$} \geq n \cdot  const(S)  \mbox { with } const(S)>0  $$
if  $S$ is the $m$-sphere $S^m(1)$,  $m\geq 2$, or, more generally, if  the   iterated loop space  
$\Omega^k(S)$ for some $k\geq 1$  has non-zero rational homology groups  $H_i(\Omega^k(S);\mathbb Q)$ for $i$ from a subset of positive density in $\mathbb Z_+$.
\vspace {1mm}

(ii) {\it Lipschitz maps.} Let  $X$ and $S$ be metric spaces and $Lip_\lambda(X$\hspace {-0.5mm}$\to$$S)$ be the space of $\lambda$-Lipschitz maps with the above length metric which we now denote $dist_\lambda$ and where we observe that
the   inclusions  
$$\mbox {$Lip_{\lambda_1}((X$\hspace {-0.5mm}$\to$$S), dist_{\lambda_1})$ $\subset$ $Lip_{\lambda_2}((X$\hspace {-0.5mm}$\to$$S),dist_{\lambda_2})$}, \hspace {1mm} \lambda_1 
\leq \lambda_2,$$
are distance decreasing.

The simplest space here, as earlier is where $X$ is a ball, but now the geometry of the ball  is essential. For instance, 
$$diam_\lambda\mbox {($Lip_{\lambda}(B^n(R)$\hspace {-0.5mm}$\to$$S))$}\leq \lambda R+diam(S),$$
where $B^n(R)$ is the  Euclidean $R$-ball, where $diam_\lambda$ is the diameter measured with $dist_\lambda$ and where, observe, 
$$\mbox {$Lip_{\lambda_1}(B^n(R_1)$\hspace {-0.5mm}$\to$$S) $ = $Lip_{\lambda_2}(B^n(R_2)$\hspace {-0.5mm}$\to$$S)$} \hspace {1mm} \mbox { for    $\lambda_1R_1=\lambda_2 R_2$}.$$
 
More interesting is the lower bound
$$diam_{c\lambda}\mbox {($Lip_{\lambda}(B^n(R)$\hspace {-0.5mm}$\to$$S))$}\geq const (S, c)\lambda R,$$
 which holds whenever {\it the real homology $H_i(S,\mathbb R)$ does not vanish for some $i\leq n$.}

In fact, 
$$diam_{c\lambda}\mbox {($Lip_{\lambda}(B^n(R)$\hspace {-0.5mm}$\to$$S))$}\geq diam_{c\lambda}\mbox {($Lip_{\lambda}(B^i(R)$\hspace {-0.5mm}$\to$$S))$  for $n\geq i$},$$
while evaluation  $f\mapsto h(f)$ of a  non-cohomologous to zero real $i$-cocycle  $h$
at  $1$-Lipschitz maps $f:B^i(R)\to S$ \footnote{Such an  $h$ in our case is a differential $i$-form on $S$ and its "evaluation"
is the integral of this form over $B^i$ mapped to $S$ by $f$.}  defines a  {\it $C\cdot R^{n-1}$-Lipschitz map} from $Lip_{1}(B^i(R)$\hspace {-0.5mm}$\to$$S))$ to $\mathbb R$  for some $C=C(S, c) )$. It follows that the 1-Lipschitz maps $f$ with $h(f) \approx  vol_i (B^i(R)) \approx R^n$  --  these  exist by the Hurewicz-Serre theorem  for the {\it minimal} $i$ where $H_i(S,\mathbb R)\neq 0$  -- are within distance $\gtrsim R$ from the constant maps.\vspace {1mm}

\hspace {40mm}{\large \it  Three  Questions.}\vspace {1mm}

    [a] {\sf Are the diameters  $$\mbox {$diam_c(Lip_{1}(B^n(R)$\hspace {-0.5mm}$\to$$S)))$}$$ bounded  for a large  fixed $c$ and  $R\to \infty$   if  $H_i(S,\mathbb R)=0$ for $i=1,2,...,n$.}

(It is not hard to show that $diam_{1}(Lip_{\lambda}(B^1(R)$\hspace {-0.5mm}$\to$$S^m(1)))\leq const\cdot \log(R)$)
\vspace {1mm}

  [b] {\sf  What is the   asymptotics  of  the diameters  $$\mbox {$diam_c(Lip_{1}(B_H^n(R)$\hspace {-0.5mm}$\to$$S)))$}$$ for the  hyperbolic balls $B_H^n(R)$ and   $R\to \infty$}?\vspace {1mm}

 [c] Let $S$ be a Riemannian manifold homeomorphic to  the connected sum of  twenty copies of $S^2\times S^2$.

{\sf Are there 1-Lipschitz maps $f_R : B^4(R) \to S$, $R\to \infty$, such that
 $h(f_R)\geq const\cdot R^4$ for a   cocycle $h$   (e.g. a closed  4-form) which represents the fundamental cohomology class $[S]\in H^4(S; \mathbb R)$, and some  $const=const(S)>0$?}\vspace {1mm}

\section {Thin Mean Convex}

Start by observing that 

{\sl smooth domains $V\subset \mathbb R^n$ with $mean.curv(\partial V)\geq 0$ are diffeotopic to  regular neighbourhoods of  subpolyhedra $P\subset U$ with $codim (P)\geq 2$.}\vspace {1mm}

In fact, if such a $V$ is bounded, this follows 
 by Morse theory applied to a linear function on $V$, while an unbounded  $V$ can be exhausted by bounded $V_i\subset Y$  with $min.curv(\partial V_i)>0$.\vspace {1mm}

Conversely, smooth submanifolds and, more general, piecewise smooth polyhedral subsets of codimension $\geq 2$ in $\mathbb R^n$ possess arbitrary thin mean convex  neighbourhoods.

In fact, the "staircase" surgery construction  for manifolds with positive scalar curvature
 (see \cite  {GL1981},  \cite {BDS2018})   applied to mean curvature  allows an attachment of such thin domains to thick ones, as follows. \footnote{It is worth remembering that the natural continuous Riemannian metrics on the doubles of a  smooth domain $V\subset \mathbb R^n$ with 
 $mean.curv (\partial V) >0$ admits  arbitrarily fine approximations by smooth metrics with $Sc>0$. Thus all shapes  and constructions  you encounter with   $mean.curv>0$ are  also seen with   $Sc>0$.}
 \vspace {1mm}

 \textbf{ [$\CIRCLE$\hspace {-0.7mm}{\large $\Yleft$}\hspace {-0.6mm}$_\Yleft^\Yleft$]}. {\sf Let $X$ be a Riemannin manifold, $V\subset X$ a smooth domain and   $P\subset X$  a piecewise smooth polyhedral  subset. 
Let $\phi(x)$, $x\in X$, be a continuous function such that  

$\bullet$  \hspace{1mm}  $\phi(x)\leq mean.curv(\partial V,x)$ for all $x\in \partial X$;

$\bullet$  \hspace{1mm}  $\phi(x)<mean.curv(\partial V,x)$ for all $x\in \partial X\cap P$.}

{\sl Given a   neighbourhood  $U\subset X$ of 
the closure of the difference $P\setminus V$,
 there exists a smooth domain $V'$ in $X$, such that 

$\ast$  \hspace{1mm} $U\supset V' \supset V\cup P;$ 

$\ast$  \hspace{1mm} $V'\setminus U=V\setminus U;$ 

$\ast$  \hspace{1mm}  $mean.curv(\partial V',x)\geq \phi(x)$for all $x\in \partial V'$.

Moreover,  if $P$ is  transversal to the boundary $\partial V$ then there is such a  $ V'$, whose intersection with $U$  serves as a {\sf regular neighbourhood} of $V\cap P$ intersected with $U$. In particular, $V'$ homotopy retracts to $V\cup P$.}

 \vspace {1mm}


This shows that, unlike to what happens to "thick" mean convex  domains, there  are  few (if at all)  {\it global} restrictions on  the shapes of the  "thin" ones 
 but there are  plenty of local ones,  where an essential  point is to understand which domains should be qualified   as "very thin". Below are some definitions and observations which may   clarify this point.

\vspace {1mm}

 \textbf {Thin$\mathbf{_{mean>0}}$}. Given a closed  subset $Y$ in a Riemannian $n$-manifold $X$, define 
 $vol_\partial (Y)$  as the infimum of the $(n-1)$-volumes of the  boundaries  of arbitrarily small  neighbourhoods $U\supset Y$  of $Y$ in $X$, i.e. 
 $$vol_\partial (Y)=\liminf_{ U\to Y} vol_{n-1}(\partial U).$$ 
Thus, \vspace {1mm}

{\sl $vol_\partial (Y)<\alpha$ if and only there exist arbitrary small neighbourhoods $U\supset Y$

 with
 $vol_{n-1}(\partial U)<\alpha$.}\vspace {1mm}

Observe that the so defined  $vol_{\partial}$ is bounded by the $(n-1)$-dimensional  Hausdorff measure,
$$ vol_\partial (Y)\leq mes_{n-1}(Y).$$
In particular,\vspace {1mm}

{\it closed subsets $Y\subset X$  with vanishing $(n-1)$-dimensional  Hausdorff measure  

 have  $vol_\partial (Y)=0$}
  (but the converse, probably, is not true). \vspace {1mm}
  
 Next, write $$mean.curv_\partial(Y) \geq \mu_0$$
 if for all $\varepsilon>0$ all  neighbourhoods  $U'\supset Y$ contain smaller smooth neighbourhoods $U\supset Y$
 such that 
 $$mean.curv (\partial U)\geq \mu_0 -\varepsilon$$ 
 
 \textbf {Implication $\mathbf {vol_\partial (Y)=0\Rightarrow mean.curv_\partial(Y)=\infty}.$} To show this let $\mu=\mu(x) $ be a  continuous function  on $U'\setminus Y$  which is  $\geq \mu_0$ for a given $\mu_0$  and which may  blow  up at $Y$.  
 
 Let $U_0$ be a  {\it $\mu-bubble$}
pinched between $Y$ and $U'$, i.e.  $U_0$ minimises the following functional
$$U\mapsto vol_{n-1}(\partial U)-\int_U\mu(x)dx.$$
 
 If $Y$ is compact and $ \mu$ is sufficiently large near $X$  such a $\mu$-bubble $U_0$ exists and
 
 $\bullet_1$ the boundary $\partial U_0$ is
  smooth away from a possible singular  subset  $\Sigma\subset \partial U_0$ of codimension $\geq 7$;
 
 $\bullet_2$    $mean.curv(\partial U_0,x)=\mu(x)$ at the regular  points $x\in \partial U_0$;
  
   $\bullet_3$ $U_0$ can be approximated by domains $U$ with {\it smooth} boundaries $\partial U$  such that $mean.curv (\partial U)\geq \mu-\varepsilon$ for a given $\varepsilon>0$. 
   
(These $\bullet_1$ and   $\bullet_2$ are standard results of the geometric measure theory and 
 $\bullet_3$ is an elementary exercise, see \cite  {Plateau-Stein}  and  \cite  {billiards} for  details.\footnote {I apologise for referring to my own articles, but I could not find  what is needed on the web except for  a 1987-paper  by F.H. Lin:
  {\sl Approximation by smooth embedded hypersurfaces with positive mean curvature}, Bull.
 Austral.Math.Soc 36 (1987) 197-208, where only a special case is treated.}

 Probably, the implication  $\mathbf {vol_\partial (Y)=0\Rightarrow mean.curv_\partial(Y)=\infty}$ remains valid for all closed  subsets in  $X$, but the above argument, as it stands, delivers   the following  weaker property  in the non-compact case.




 \vspace {1mm}

{\it If $X$  has  uniformly bounded geometry\footnote{ This means that there exist   $\rho>0$  and $\lambda>0$, such that all  $\rho$-balls  $B_x(\rho)\subset X$ 
are $\lambda$-bi-Lipschitz homeomorphic to a Euclidean ball.}   
   then every  closed  subsets   $Y\subset X$ with $vol_\partial (Y)=0$   is  equal to  the intersection of a decreasing family of     domains  $U_\mu\subset X$, $\mu\to \infty$, where
   $mean.curv (\partial U_\mu)\geq \mu  $.}

\vspace {1mm}

{\it Remark.}  The  role  of bounded geometry is to ensure  a lower bound  on the volumes of  balls in 
the $\mu$-bubble away from $Y$  where $\mu$ is small  and, thus,   keep domains $U$ which minimise the function   $U\mapsto vol_{n-1}(\partial U)-\int_U\mu(x)dx$   within an $\varepsilon$-neighbourhood of   $Y$.

\vspace {1mm} 

 {\it Exercises.}  (a)  Let $X$ be a Riemannian manifold isometrically acted upon by a group $G$ and let  $P\subset X$ be a $G$-invariant  piecewise smooth polyhedral subset of codimension $2$.  

  Show that $P$ admits  a $G$-invariant  arbitrarily  small {\it regular} neighbourhoods with arbitrarily large mean curvatures of their boundaries.
 
 \vspace{1mm}
 
 (b)  Show that every bounded smooth domain $U_0\subset \mathbb R^n$, $n\geq 3$  contains arbitrarily long simple curves (arcs) $C\subset U$, such that
 $ curvature(C)\leq const=const(U)$.
 
Construct   such curves with  regular neighbourhoods $U\subset U_0$  which  fill almost   all of $U_0$  and such that the mean curvatures of the boundaries  $
\partial U$  tend to infinity.

Apply this successively to  $C_i\subset U_{i-1}$ and obtain
 $$U_0\supset U_1\supset  ... \supset  U_{i-1}  \supset  U_i \supset ... .$$
 such  that 
  
  \hspace {20mm}  $mean.curv(\partial U_i)\to \infty$   and   $vol(U_i)\geq const>0.$ 
  
  
  Show that the intersection $Y_\infty=\cap_iU_i$  is a compact set such that 
  
$\bullet$ $Y$  has  positive  Lebesgue measure, $vol_n(Y)>0$;
  
  $\bullet$  $vol_\partial(Y_\infty)=\infty$;
    
  $\bullet$  $ mean.curv_\partial (Y)=\infty$;
    
 $\bullet$     the topological dimension  of $Y$ is one.\vspace {1mm}

 (c) Construct similar $Y\subset \mathbb R^n$ with  $dim_{top}=m$ for all $m\leq n-2$.

\vspace {1mm}


(d) Show that the function  $\lambda |\sin  t |$, can be uniformly approximated,   for all $\lambda>0$, by 
  $C^2$-functions 
 $\varphi(t)>0$, such that the hypersurface $H_\varphi\subset \mathbb R^{n+1}$, $n\geq 2$, obtained 
 by rotating the graphs of $\varphi(t)$ around the $t$-axis in $\mathbb R^{n+1}$  has $mea.curv(H_\varphi)>0$ and, moreover,  $Sc(H_\varphi)>0$ if $n\geq 3$.

(e) Construct  Cantor (compact  0-dimensional) subsets  $Y$ in the plane which {\it are not } intersections of locally convex subsets, i.e. disjoint unions of convex ones.
\vspace {1mm}

(f)  Define "random"  Cantor  sets in $ \mathbb R^n$ of positive measure, and show for $n\geq 2$ that they  are not intersections of smooth mean convex  domains. Then, do the same for Cantor sets in $\mathbb R^n$  with the Hausdorff dimensions  $>n-1$.
  
 {\it Admission.}  Frankly, I am not 100\% certain as I haven't seriously tried to solve this exercise.

\end{document}